\documentclass[10pt]{article}
\usepackage{amsmath,amssymb,amsfonts}
\newtheorem{theorem}{Theorem}
\newtheorem{lemma}{Lemma}

\date{}

\title{On the Distribution of the Number of Goldbach Partitions of a Randomly
Chosen Positive Even Integer}
\bigskip

\author{{\bf Ljuben Mutafchiev}\\
American University in Bulgaria, 2700 Blagoevgrad, Bulgaria \\ and
Institute of Mathematics and Informatics of the \\ Bulgarian
Academy of Sciences
\\ \tt {ljuben@aubg.bg}}

\begin{document}
\maketitle

\begin{abstract}
Let $\mathcal{P}=\{p_1,p_2,...\}$ be the set of all odd primes
arranged in increasing order. A Goldbach partition of the even
integer $2k>4$ is a way of writing it as a sum of two primes from
$\mathcal{P}$ without regard to order. Let $Q(2k)$ be the number
of all Goldbach partitions of the number $2k$. Assume that $2k$ is
selected uniformly at random from the interval $(4,2n], n>2$, and
let $Y_n=Q(2k)$ with probability $1/(n-2)$. We prove that the
random variable $\frac{Y_n}{n/\left(\frac{1}{2}\log{n}\right)^2}$
converges weakly, as $n\to\infty$, to a uniformly distributed
random variable in the interval $(0,1)$. The method of proof uses
size-biasing and the Laplace transform continuity theorem.
\end{abstract}

\vspace{.5cm}



\vspace{.2cm}

\section{Introduction}

Let $\mathcal{P}=\{p_1,p_2,...\}$ be the sequence of all odd
primes arranged in increasing order. A Goldbach partition of the
even integer $2k>4$ is a way of writing it as a sum of two primes
$p_i,p_j\in\mathcal{P}$ without regard to order. The even integer
$2k=p_i+p_j, i\ge j,$ is called a Goldbach number. Let $Q(2k)$
denote the number of the Goldbach partitions of the number $2k$.
In 1742 C. Goldbach conjectured that $Q(2k)\ge 1$ for all $k>2$.
This problem still remains unsolved (for more details, see, e.g.,
[4; Section 2.8, p. 594], [8; Section 4.6] and [9; Chapter VI]).
Let $\Sigma_{2n}$ be the set of all Goldbach partitions of the
even integers from the interval $(4,2n], n>2$. The cardinality of
this set is obviously
\begin{equation}\label{sigma}
\mid\Sigma_{2n}\mid=\sum_{2<k\le n} Q(2k).
\end{equation}

In this paper, we consider two random experiments. In the first
one, we select a partition uniformly at random from the set
$\Sigma_{2n}$, i.e., we assign the probability
$1/\mid\Sigma_{2n}\mid$ to each Goldbach partition of an even
integer from the interval $(4,2n]$. An important statistic (random
variable) of this experiment is the Goldbach number $2G_n\in
(4,2n]$ partitioned by this random selection. Its probability mass
function (pmf) is given by
\begin{equation}\label{gpmf}
f_{G_n}(x)=\frac{Q(2k)}{\mid\Sigma_{2n}\mid} \quad \text {if
$x=2k\in (4,2n]$},
\end{equation}
and zero elsewhere.

{\it Remark 1.} It was established in [7] that
\begin{equation}\label{equivalence}
\mid\Sigma_{2n}\mid\sim\frac{2n^2}{\log^2{n}}, \quad n\to\infty.
\end{equation}
The proof is based on a classical Tauberian theorem due to
Hardy-Littlewood-Karamata [3; Chapter 7]. Recently, in a private
communication, Kaisa Matom\"{a}ki [5] showed me a shorter and
direct proof of (\ref{equivalence}) that uses only the Prime
Number Theorem [4; Section 17.7] and partial summation.

In the second random experiment, we select an even number $2k\in
(4,2n]$ with probability $1/(n-2)$. Let $Y_n$ be the statistic,
equal to the number $Q(2k)$ of its Goldbach partitions. Obviously,
the pmf of $Y_n$ is
\begin{equation}\label{ypmf}
f_{Y_n}(x)=\frac{1}{n-2} \quad \text {if $x=Q(2k), \quad 2k \in
(4,2n]$},
\end{equation}
and zero elsewhere.

The main goal of this paper to study the limiting distribution of
the random variable $Y_n$. In Section 2 we show that $Y_n$,
appropriately normalized, converges weakly, as $n\to\infty$, to a
random variable that is uniformly distributed in the interval
$(0,1)$.

The method of proof uses size biasing and the Laplace transform
continuity limit theorem [2; Chapter XIII, Section 1].

\section{The Limiting Distribution of $Y_n$}

The first step is to determine the asymptotic of the expected
value of $Y_n$. From (\ref{sigma}), (\ref{equivalence}) and
(\ref{ypmf}) it follows that
\begin{equation}\label{ev}
\mathbb{E}(Y_n) =\frac{1}{n-2}\sum_{2<k\le n}
Q(2k)\sim\frac{2n}{\log^2{n}}=\frac{1}{2} b_n, \quad n\to\infty,
\end{equation}
where
\begin{equation}\label{bn}
b_n =\frac{n}{\left(\frac{1}{2}\log{n}\right)^2}, \quad n>2.
\end{equation}

Next, we will introduce the concept of size-biasing. The
definition we present below is given in [1; Section 4.2]. Suppose
that $X$ is a non-negative random variable with finite mean $\mu$
and distribution function $F$. The notation $X^*$ is used to
denote a random variable with distribution function given by
\begin{equation}\label{sb}
F^*(dx) =\frac{xF(dx)}{\mu}, \quad x>0.
\end{equation}
The random variable $X^*$ and the distribution function $F^*$ are
called size-biased versions of $X$ and $F$, respectively.

Our goal is to show that the size-biased version $Y_n^*$ of $Y_n$
is $G_n$. To see this, we set in the right-hand side of (\ref{sb})
$x=Q(2k)$ $(2k\in(4,2n])$ and $\mu=\mathbb{E}(Y_n)$. Using
(\ref{ev}), (\ref{ypmf}), (\ref{gpmf}) and (\ref{sigma}), we
obtain
\begin{equation}\label{sbpmf}
  f_{Y_n^*}(2k)=\frac{\frac{1}{n-2} Q(2k)}{\mathbb{E}(Y_n)}=\frac{\frac{1}{n-2} Q(2k)}
{\frac{1}{n-2} \sum_{2<k\le n}Q(2k)}
=\frac{Q(2k)}{\mid\Sigma_{2n}\mid} =f_{G_n}(2k).
\end{equation}
If we multiply both sides of (\ref{sbpmf}) by $e^{-2\lambda
k/(2n)}=e^{-\lambda k/n}, \lambda>0,$ and sum up over $k\in
(2,n]$, we observe that
$$
\sum_{2<k\le n}\frac{e^{2\lambda k/(2n)}\left(\frac{1}{n-2}\right)
Q(2k)}{\mathbb{E}(Y_n)} =\frac{\mathbb{E}(Y_n e^{-\lambda X_n})}
{\mathbb{E}(Y_n)} =\mathbb{E}(e^{-\lambda G_n/n}),
$$
where $X_n$ denotes a random variable that assumes the values
$2k/(2n)=k/n, k\in (2,n]$, with probability $1/(n-2)$. Hence, for
fixed $n>2$, we have
$$
\mathbb{E}(Y_n e^{-\lambda
X_n})=(\mathbb{E}(Y_n))\mathbb{E}(e^{-\lambda G_n/n}).
$$
Let
\begin{equation}\label{zn}
Z_n=\frac{1}{b_n}Y_n, \quad n>2,
\end{equation}
where the scaling factor $1/b_n$ is defined by (\ref{bn}).
Clearly, $Z_n$ satisfies the same identity as $Y_n$. We have
 \begin{equation}\label{sbexpect}
\mathbb{E}(Z_n e^{-\lambda
X_n})=(\mathbb{E}(Z_n))\mathbb{E}(e^{-\lambda G_n/n}), \quad n>2.
\end{equation}

We start our asymptotic analysis of (\ref{sbexpect}) with
\begin{equation}\label{ezn}
\lim_{n\to\infty}\mathbb{E}(Z_n)=\frac{1}{2},
\end{equation}
which follows from (\ref{ev}) and (\ref{zn}). The second factor in
the right-hand side of (\ref{sbexpect}) is the Laplace transform
of the random variable $G_n/n$. Its limiting distribution, as
$n\to\infty$, was found in [7]. (The proof is based on the
asymptotic equivalence (\ref{equivalence}).) We state this result
in the following separate lemma.

\begin{lemma}
The sequence of random variables $\{G_n/n\}_{n\ge 2}$ converges
weakly, as $n\to\infty$, to the random variable
$U=\max{\{U_1,U_2\}}$, where $U_1$ and $U_2$ are two independent
copies of a uniformly distributed random variable in the interval
$(0,1)$.
\end{lemma}

{\it Remark 2.} Clearly, the probability density function of the
random variable $U$ equals $2x$, if $0<x<1$, and zero elsewhere.
The $r$th moment of $U$ is $\frac{2}{r+2}, r=1,2,...$. The Laplace
transform of the random variable $U_1$ (uniformly distributed in
the interval $(0,1)$) is
\begin{equation}\label{laplaceuone}
\varphi(\lambda) =\frac{1-e^{-\lambda}}{\lambda}, \quad \lambda>0,
\end{equation}
while the Laplace transform of $U$ is
$$
\frac{2}{\lambda^2}(1-e^{-\lambda}-\lambda e^{-\lambda})
=-2\varphi^\prime(\lambda), \quad \lambda>0.
$$

Further, for any fixed $s>0$, we integrate (\ref{sbexpect}) with
respect to $\lambda$ over the interval $(0,s)$. Applying Fubini's
theorem [2; Chapter IV, Section 3], we obtain
\begin{equation}\label{integral}
\mathbb{E}\left(\frac{Z_n}{X_n}(1-e^{-sX_n})\right)
=(\mathbb{E}(Z_n))
\mathbb{E}\left(\frac{1-e^{-sG_n/n}}{G_n/n}\right).
\end{equation}
To find the limit of the right-hand side of (\ref{integral}), we
combine the result of Lemma 1 with the continuity theorem for
Laplace transforms (see, e.g., [2; Chapter XIII, Section 1]) and
the Lebesgue dominated convergence theorem. Using the probability
density function of the random variable $U$, (\ref{ezn}) and
(\ref{laplaceuone}), we deduce that
\begin{equation}\label{rhs}
\lim_{n\to\infty}(\mathbb{E}(Z_n))
\mathbb{E}\left(\frac{1-e^{-sG_n/n}}{G_n/n}\right)
=\frac{1}{2}\mathbb{E}\left(\frac{1-e^{-sU}}{U}\right)
=1-\varphi(s).
\end{equation}
The random variable $X_n$ in the left-hand side of
(\ref{integral}) has simple probabilistic interpretation: $nX_n$
equals an integer, chosen uniformly at random from the interval
$(2,n]$. So, it converges weakly, as $n\to\infty$, to the random
variable $U_1$ (uniformly distributed in the interval $(0,1)$).
Applying again the continuity theorem for Laplace transforms [2;
Chapter XIII, Section 1], we have
\begin{equation}\label{phi}
\lim_{n\to\infty}\mathbb{E}(e^{-sX_n})=\varphi(s), \quad s>0,
\end{equation}
where $\varphi$ is given by (\ref{laplaceuone}). Moreover,
(\ref{integral}) and (\ref{rhs}) imply that
\begin{equation}\label{lhs}
\lim_{n\to\infty}\mathbb{E}\left(\frac{Z_n}{X_n}(1-e^{-sX_n})\right)
=1-\varphi(s).
\end{equation}
Now, it is not difficult to show that $Z_n/X_n\to 1$ in
probability, as $n\to\infty$. In fact, if this is not true, then
there is a number $\epsilon\in (0,1)$, such that, for infinitely
many values of $n$, either $Z_n/X_n\le 1-\epsilon$ or $Z_n/X_n\ge
1+\epsilon$ for these values of $n$ with probability tending to
$1$. Both inequalities contradict with (\ref{phi}) and (\ref{lhs})
for these $n$. Hence, for any $\eta>0$,
$$
\lim_{n\to\infty}\mathbb{P}\left(\mid\frac{Z_n}{X_n}-1\mid\ge
\eta\right)=0,
$$
which implies that $Z_n$ and $X_n$ have one and the same limiting
distribution. Thus we obtain the following limit theorem.

\begin{theorem} The sequence $\{Z_n=Y_n/b_n\}_{n>2}$, with $b_n$
given by (\ref{bn}), converges weakly, as $n\to\infty$, to the
random variable $U_1$, which is uniformly distributed in the
interval $(0,1)$.
\end{theorem}

Theorem 1 shows that the number of Goldbach partitions of even
integers $2k\in (4,2n]$ is typically of order
$$
ab_n=\frac{4an}{\log^2{n}},
$$
where $0<a<1$ (see (\ref{bn})). This informal evidence in favor of
the Goldbach's conjecture does not give us any rigorous argument
to prove it. Finally, we note that in number theory special
interest is paid on asymptotic estimates for the probability
$\mathbb{P}(Y_n=0)$ (see, e.g., [9; Chapter VI]). For instance,
Montgomery and Vaughan [6] have shown that, for sufficiently large
$n$, there exist a positive (effectively computable) constant
$\delta<1$, such that $\mathbb{P}(Y_n=0)\le (2n)^{-\delta}$.

\end{document}